\documentclass[11pt, a4paper, reqno]{amsart}
\usepackage{latexsym}
\usepackage{amsthm, amssymb, amsfonts}
\usepackage{amsmath}
\usepackage[colorlinks=true,linkcolor=blue,citecolor=dgreen]{hyperref}
\usepackage[shortlabels]{enumitem}
\usepackage{amssymb}
\usepackage{amsmath,float}
\usepackage{tikz}
\usepackage{cleveref}
\usepackage[numbers]{natbib}
\usepackage[english]{babel} 
\usepackage{blindtext}

\DeclareMathOperator{\col}{col}
\definecolor{shadecolor}{rgb}{0.9,0.9,0.9}
\def\={\equiv}

\newcommand{\myblueshadedcell}[1]{\tikz[baseline,remember picture]{\node[fill=blue!60,anchor=base] (cell) {#1};}}
\newcommand{\myredshadedcell}[1]{\tikz[baseline,remember picture]{\node[fill=red!60,anchor=base] (cell) {#1};}}

\newtheorem{theorem}{Theorem}[section]

\newtheorem{lemma}[theorem]{Lemma}

\theoremstyle{definition}

\newcommand{\Z}{\mathbb{Z}}
\renewcommand{\AA}{\mathcal{A}}

\newcommand{\HH}{\mathcal{I}}
\newcommand{\LL}{\mathcal{L}}

\newcommand{\TT}{\mathcal{T}}
\newcommand{\UU}{\mathcal{U}}
\newcommand{\VV}{\mathcal{V}}

\definecolor{dgreen}{rgb}{0,.8,0}
\renewcommand{\leq}{\leqslant}
\renewcommand{\geq}{\geqslant}
\renewcommand{\le}{\leqslant}
\renewcommand{\ge}{\geqslant}
\renewcommand{\emptyset}{\varnothing}

\def\eref#1{$(\ref{#1})$}

\def\Tref#1{Table~$\ref{#1}$}

\crefname{equation}{}{}
\crefname{enumi}{}{}
\crefmultiformat{lemma}{lemmas~#2#1#3}{ and~#2#1#3}{, #2#1#3}{ and~#2#1#3}
\title{\bf Latin Squares whose transversals share many entries}

\author{Afsane Ghafari and Ian M. Wanless}

\thanks{
School of Mathematics, Monash University, Vic 3800, Australia.
{\tt afsane.ghafaribaghestani@monash.edu, ian.wanless@monash.edu}.
}

\date{}
\begin{document}

\maketitle

\begin{abstract}
We prove that, for all even $n\geq10$, there exists a latin square of
order $n$ with at least one transversal, yet all transversals coincide
on $ \big\lfloor n/6 \big\rfloor$ entries. These latin squares have at
least $ 19 n^2/36 + O(n)$ transversal-free entries. We also prove that
for all odd $m\geq 3$, there exists a latin square of order $n=3m$
divided into nine $m\times m$ subsquares, where every transversal hits
each of these subsquares at least once.
\end{abstract}

\section{Introduction}

A \emph{latin square} is an $n \times n$ array consisting of $n$
distinct symbols where each symbol appears exactly once in each row
and column.  For any set $\HH$ of cardinality $n$, the rows and
columns of a latin square $L$ can be indexed by $\HH$ so that $L$ can
be defined by a collection of $n^2$ ordered triples
$\{(r,c,s)\}\subseteq\HH(L)^3$, where $r$ and $c$ denote,
respectively, the row and the column indices of a cell containing the
symbol $s$. We say that $L$ is \emph{indexed} by $\HH(L)$. In this
representation $(r,c,s)$ is considered an \emph{entry} in $L$, and the
latin property ensures that two different entries agree in at most one
coordinate.

A \emph{transversal} of a latin square $L$ of order $n$ is a set of
$n$ entries of $L$ in which each row, column and symbol is represented
exactly once.  We say that an entry of a latin square is \emph{pinned} if
it is contained in every transversal (and there is at least one
transversal).
We refer to an entry that is not in any transversal as a
\emph{transversal-free} entry.  We also use the notation $\tau(L)$ for
the number of transversal-free entries in $L$. We say two transversals
are \emph{disjoint} if they share no entries.

Our main results in this article are the following:

\begin{theorem}\label{verymaintheorem}
  Let $n\geq 10$ be an even integer. There exists a latin square of
  order $n$ with at least $\big\lfloor {n}/{6} \big\rfloor$ pinned
  entries.
\end{theorem} 

\begin{theorem}\label{resultfrommain}
  For all even $n \geq 10$, there exists a latin square $L$ of order $n$
  satisfying $n^2>\tau(L) \geq 19 n^2/36  + O(n)$.
\end{theorem}

\begin{theorem}\label{mySecTheorem}
  For all odd $m\geq 3$, there exists a latin square of order $n=3m$ composed of nine subsquares of order $m$, where every transversal hits each of these subsquares at least once.
\end{theorem}

We give constructive proofs of \Cref{verymaintheorem,resultfrommain}
in \Cref{FirstSec} and \Cref{mySecTheorem} in \Cref{secsec}.  Egan and
Wanless \citep{egan2012latin} proved a weaker version of
\Cref{verymaintheorem}. Namely, for all even $n \geq 10$, there exists
a latin square of order $n$ containing at least one pinned entry.  To
our knowledge, the lower bound given in \Cref{resultfrommain} is the
best bound among latin squares of even order which include at least
one transversal. Egan and Wanless \citep{egan2012latin}
observed that for all $n \geq 4$, there exists a latin square $L$ of
order $n$ with $\tau(L) \geq 7$. They also proved the following theorem:

\begin{theorem}\label{t:quadodd}
For all odd $m \geq 3$, there exists a latin square of order $3m$ that
contains an $(m-1) \times m$ latin subrectangle consisting of
transversal-free entries.
\end{theorem}

The latin squares in \Cref{t:quadodd} represent the only previously
known family that have transversals, but have more than a constant
fraction of their entries not in any transversal.  It is important to
note, however, that there exist latin squares of even order that
contain no transversal. One of the earliest theorems regarding
transversals was proven by Euler \citep{euler1782recherches}, who
showed that the addition table of $\Z_n$ has no transversals when $n$
is even.  In fact there are a huge number of latin squares of even
order that have no transversals.  Cavenagh and Wanless
\citep{cavenagh2017latin} showed that for even $n \rightarrow \infty$,
there are at least $n^{n^{3/2} (1/2-o(1))}$ species of
transversal-free latin squares of order $n$.  Here, a \emph{species}
(also called \emph{main class}) is an equivalence class of latin
squares under natural operations that preserve much of the
combinatorial structure, including the number of transversals.  Note
that $\tau(L)=n^2$ for latin squares that have no transversals.

\section{Transversals coinciding on several entries}\label{FirstSec}

In this section we give a constructive proof of
\Cref{verymaintheorem,resultfrommain}.  First we introduce the
\emph{$\Delta$-Lemma} which first appeared independently
in \citep{evans2006latin} and \citep{WW06}.
Consider the function $\Delta : L \rightarrow \Z$ to be the following:
\begin{align}\label{Copyfascinating}
\Delta (r,c,s) \equiv s-r-c \bmod n\text{ where } -\dfrac{n}{2} < \Delta(r,c,s) \leq \dfrac{n}{2}.
\end{align}

\begin{lemma}\label{l:Delta}
Let $L$ be a latin square of order n indexed by $\Z_n$. If $T$ is a transversal of $L$ then, modulo $n$,
$$\sum_{(r,c,s)\in T} \Delta(r,c,s) \equiv
\begin{cases}
	0  & \text{if $n$ is odd},\\
	\frac{1}{2} n & \text{if $n$ is even}.
\end{cases}$$
\end{lemma}
Our latin squares are indexed
by a cyclic group $\Z_n$, and also have their symbols
chosen from $\Z_n$. All calculations are in $\Z_n$
unless otherwise stated and we assume our residues to be
$\{0,1,\ldots, n-1\}$.

For a latin square $L$, the notation $\col(a)$ is used to specify
that transversal $T$ includes the entry $(a,\col(a),L[a,\col(a)])$.

We define a \emph{suitable diagonal} of a latin square $L$ of even
order $n$ to be any set of $n$ entries from different rows and
different columns, such that the sum of the $\Delta$ values of the
chosen entries is $n/2\bmod n$.  In light of \Cref{l:Delta}, every
transversal is a suitable diagonal.  So to prove
\Cref{verymaintheorem} it suffices to prove the following slightly
stronger statement instead:

\begin{theorem}\label{t:surrogate}
  For all even $n\ge10$, there exists a latin square of order $n$ that
  contains a transversal, but in which there are $\big\lfloor n/6\big\rfloor$
  entries that occur in every suitable diagonal.
\end{theorem}

We prove \Cref{t:surrogate} by separating into three cases depending on the value of $n\bmod6$.
We start by assuming that $n= 6k$, for an integer $k\geq 2$. Let $\TT_n$
be the latin square of order $n$ defined by,
\begin{equation}\label{Structure}
  \TT_n[a,b] = \begin{cases}
    a+b+2 &  \text{if  $(a,b) \in \{(0,2),(1,0)\},$} \\  
    a+b+1& \text{if $(a,b) \in \{(0,1),(2,1)\},$}\\
    a+b-1 & 	\text{if $(a,b)\in \{(1,1), (1,2), (2,2), (3,1)\},$}\\
    a+b-2  &   \text{if $(a,b) = (3,0),$} \\
    a+b+3  &   \text{if $b> 1 ,  b\equiv 1$  mod 3  and  $a=0,$}  \\
    a+b-3  &   \text{if $b> 1 , b\equiv 1$  mod 3 and $a=3,$ }  \\
    a+b-2  &   \text{if $4\le a\le3k-3$, $a\equiv 0$ mod 3 and $b\equiv 0$ mod 2,} \\ 
    a+b+2  & \text{if $4\le a\le3k-3$, $a\equiv 1$ mod 3 and $b\equiv 0$  mod 2, $b\neq n - 2a + 2,$} \\ 
    a+b+1  &   \text{if $4\le a\le3k-3$, $a\equiv 1$  mod 3 and $b\in \{n - 2a + 3,  n - 2a + 2\},$}  \\
    a+b+1   &  \text{if $4\le a\le3k-3$, $a\equiv 2 $  mod 3 and $b= n - 2a + 4,$} \\
    a+b-1  &   \text{if $4\le a\le3k-3$, $a\equiv 2$  mod 3 and $b=n - 2a + 5,$}  \\
    a+b & \text{otherwise.}
  \end{cases}
\end{equation}

\begin{lemma}\label{zeromodsix}
  Let $n=6k$ for $k\ge2$.
  The latin square $\TT_n$ has a transversal, and all
  suitable diagonals of $\TT_n$ include the $k$ entries
\begin{align}\label{0modd6}
(1,0,3), (2,1,4), (5,n-6,0), (8,n-12,n-3),\ldots,( 3k-4,12,3k+9).
\end{align}
\end{lemma}

\begin{proof}
From \cref{Structure} we see that $-3\le\Delta(r,c,s)\le3$ for all $(r,c,s)$ in $\TT_n$.
The sum of the minimum $\Delta$-value from each row is
$$ 0-1-1-3+\sum_{i=1}^{k-2} (0-1-2)= -3k+1 > -3k = - \dfrac{n}{2}.$$
Also, the sum of the maximum $\Delta$-value from each row is
$$3+2+1+0+\sum_{i=1}^{k-2} (2+1+0)=3k=\dfrac{n}{2}.$$
Therefore, any suitable diagonal in $\TT_n$ must include an entry with the maximum $\Delta$-value from each row, meaning that it must include the entries in \cref{0modd6}.

Now, we demonstrate the existence of a transversal $T$ in $\TT_n$, 
defined as follows:
\begin{equation*}
  \col(a) = \begin{cases}
    4  &   \text{if $a=0, $}\\
    a-1 &  \text{if $a\in \{1,2,3\}, $}\\ 
    n-2a+4   &    \text{if $4 \leq a \leq 3k-1$ and $a\not\equiv 0$  mod 3,} \\  
    n-2a+7  &    \text{if $4 \leq a \leq 3k-1$ and $a\equiv 0$ mod 3,} \\  
    n-2a+3   &    \text{if $3k\le a<n$  and $a\equiv 0$  mod 3,} \\  
    n-2a+6  &    \text{if $3k+2\le a<n$  and $a\equiv 1$ mod 3,} \\    
    n-2a+9  &    \text{if $(3k+2\le a<n$  and $a\equiv 2$  mod 3) or $a=3k+1$.}
  \end{cases}
\end{equation*}
In order to prove that $T$ is a transversal, we show that every column and symbol is represented exactly once in $T$. Let $C_{a,b}$ and $S_{a,b}$ be, respectively, the set of columns and symbols used by the entries in $T$ in the rows $\{a,a+1,\ldots,b\}$. Then
\begin{align*}
 C_{0,3}&=\{0,1,2,4\},\\
C_{4,3k-1}&=\{6s+r: 1 \leq s \leq k - 1 , 0 \leq r \leq 2\} \setminus \{7\},\\
 C_{3k,n-1}&= \{6s+r: 1 \leq s \leq k - 1 ,3 \leq r \leq 5\} \cup \{3,5,7\},
\end{align*}
verifying that $C_{0,n-1}=\{0,1,\ldots,n-1\}$.
Also,
\begin{align*}
 S_{0,3}&=\{3,4,5,7\},\\
S_{4,3k-3}&=\{s:3k+9\le s\le6k+2\}\\
S_{3k-2,3k+2}&=\{3k+3,3k+5,3k+6,3k+7,3k+8\},\\
S_{3k+3,n-1} &=  \{s:8\le s\le3k+2\}\cup\{6,3k+4\},
\end{align*}
confirming that $S_{0,n-1}=\{0,1,\ldots,n-1\}$, and completing the
proof of the lemma.
\end{proof}

Next suppose that $n=6k+2$, where $k \geq 2$ is an integer. Consider
the latin square $\UU_n$ given by:
\begin{equation}\label{SecFamily}
  \UU_n[a,b] = \begin{cases}
    a+b+1    &    \text{if $(a,b) \in \{(1,3),(2,5),(3,4)\},$}  \\ 
    a+b-1 &   \text{if $(a,b) \in \{ (1,4), (2,4), (3,5), (4,3), (4,4) \},$} \\
    a+b+2 &  \text{if $(a,b) \in \{(0,4),(2,3) \},$ } \\
    a+b+2 &  \text{if $b \neq 4, b\equiv 0$ mod 2 and  $a=2$ } \\
    a+b+3  &   \text{if $(b > 3, b\equiv 0$ mod 3 and  $a=0$)  or $(a,b)=(0,1),$ } \\
    a+b-3  &   \text{if $(b > 3 , b\equiv 0$ mod 3 and  $a=3$)   or $(a,b)=(3,1),$} \\
    a+b-2  &   \text{if $(b\neq 4,\, b\equiv 0$  mod 2 and $a = 4$)  or $(a,b)=(3,3),$ }  \\
    a+b+2  &   \text{if $5\le a\le3k-2$, $a\equiv$ 2 mod 3 and $b\equiv 0$  mod 2, $b\neq n - 2a + 4,$}  \\ 
    a+b+1  &   \text{if $5\le a\le3k-2$, $a\equiv 2$  mod 3 and $b\in \{n - 2a +4,  n - 2a + 5\},$}   \\
    a+b+1      & \text{if $5\le a\le3k-2$, $a\equiv 0$  mod 3 and $b= n - 2a + 6,$} \\
    a+b-1  &   \text{if $5\le a\le3k-2$, $a\equiv 0$  mod 3 and $b=n - 2a +7, $}  \\
    a+b-2 &	 \text{if $5\le a\le3k-2$, $a\equiv 1$  mod 3 and $b\equiv 0$ mod 2,}  \\ 
    a+b & \text{otherwise.}
  \end{cases}
\end{equation}

\begin{lemma}\label{twomodsix}
The latin square $\UU_n$, where $n=6k+2$, has a transversal, and all suitable diagonals in $\UU_n$ include the $k$ entries 
\begin{align}\label{2modd6}
(1, 3, 5), (3,4,8), (6,n-6,1), (9,n-12,n-2),\ldots,( 3k-3,14,3k+12).
\end{align}
\end{lemma}

\begin{proof}
  Based on \cref{SecFamily}, we have $-3\le\Delta(r,c,s)\le3$ for every entry $(r,c,s)$ in $\UU_n$. 
The sum of the minimum $\Delta$-value from each row is
$$ 0-1-1-3-2+\sum_{i=1}^{k-2} (0-1-2)= -3k-1 = - \dfrac{n}{2}.$$
However, the entries $(1,4,4)$ and $(2,4,5)$ both achieve the unique
minimum value of $\Delta$ within their respective rows.
No suitable diagonal can include both these entries since
they share a column.

On the other hand, the sum of the maximum $\Delta$-value from
each row is
\begin{align*}
3+1+2+1+0+\sum_{i=1}^{k-2} (2+1+0)=3k+1=\dfrac{n}{2}.
\end{align*}
Therefore, any suitable diagonal must include an entry with the
maximum $\Delta$-value from each row, which forces it to
include the entries in \cref{2modd6}.

Next, we establish the presence of a transversal $T$ in $\UU_n$,
defined by
\begin{align*}
\col(a) = \begin{cases}
  a+1  &   \text{if $a\in \{0,3,4\},$} \\
  3& \text{if $a = 1,$}\\
  8& \text{if $a = 2,$} \\
  2& \text{if $a = 3k+6,$} \\
  6& \text{if $a = 3k+1,$} \\
  7& \text{if $a = 3k+3,$}\\
  9& \text{if $a = 6k-1,$}  \\
  11& \text{if $a = 3k,$} \\
  n-2a+6   &    \text{if ($5 \leq a < 3k$ and $a\not\equiv 1$ mod 3) or $a=3k+4$,} \\  
  n-2a+9  &    \text{if $5 \leq a < 3k$ and $a\equiv 1$  mod 3,}  \\  
  n-2a+13   &    \text{if $3k+7\leq a <n$ and $a\equiv 0$  mod 3,}\\  
  n-2a+10  &    \text{if $3k+7\leq a <n$ and  $a\equiv 1$ mod 3,} \\    
  n-2a+1  &    \text{if $3k+2\leq a \leq 6k-4$  and $a\equiv 2$  mod 3.}\\  
\end{cases}
\end{align*}

In order to prove that $T$ is a transversal, we show that every column and symbol is represented exactly once in $T$.
Let $C_{a,b}$ and $S_{a,b}$ be, respectively, the set of columns and symbols used by the entries in $T$ in rows $\{a,a+1,\ldots,b\}$. Thus:
\begin{align*}
 C_{0,4}&=\{1,3,4,5,8\},\\
C_{5,3k-2}&=\{6s+r: 2 \leq s \leq k - 1 , 2\leq r\leq 4\},\\ 
C_{3k-1,3k+4} &=\{0,6,7,10,11,6k+1\},\\
C_{3k+5,n-1}&= \{6s+r: 2 \leq s \leq k-1  , r\in \{0,1,5\}\} \cup \{2,9,6k\},
\end{align*}
which together show that $C_{0,n-1} = \{0,1,\ldots,n-1\}$.
Also,
\begin{align*}
 S_{0,4}&=\{4,5,8,9,12\},\\
S_{5,3k-2}&=\{s:3k+12\le s\le6k+5\},\\
S_{3k-1,3k+4} &=\{3k+1,3k+4,3k+7,3k+9,3k+10,3k+11\},\\
S_{3k+5,n-1} &= \{3s+r:2\leq s\leq k-1,r\in\{1,5,9\}\}\cup \{6,3k+5,3k+8\},
\end{align*}
which together show that $S_{0,n-1}=\{0,1,\ldots,n-1\}$, finishing the proof of the lemma. 
\end{proof}

It remains to deal with the case $n\equiv 4 \bmod  6$.
Let $n=6k+4$, where $k \geq 1$ is an integer.
Consider the latin square $\VV_n$ given by:
\begin{equation}\label{ThirdFamily}
  \VV_n[a,b] = \begin{cases}
    a+b-1    &    \text{if $(a,b) \in \{(1,1),(1,2)\},$} \\ 
    a+b-2 &   \text{if $(a,b) \in \{(3,0), (3,2)\},$} \\
    a+b+1 &  \text{if $(a,b) \in \{(0,1)\},$}  \\  
    a+b+2 &  \text{if $(a,b) \in \{(1,0)\}, $}\\ 
    a+b+3  &   \text{if $b\equiv 2$ mod 3  and $a=0,$}   \\
    a+b-3  &   \text{if $2<b\equiv 2$ mod 3 and  $a=3,$}    \\
    a+b-2  &   \text{if $4\le a\le3k$, $a\equiv 0$  mod 3 and $b\equiv 0$  mod 2,}  \\ 
    a+b+2  &   \text{if $4\le a\le3k$, $a\equiv 1$ mod 3 and $b\equiv 0$ mod 2, $b\neq n - 2a + 2, $}  \\ 
    a+b+1   &   \text{if $4\le a\le3k$, $a\equiv 1$  mod 3 and $b\in\{ n - 2a + 2,n - 2a + 3\},$} \\
    a+b+1  &   \text{if $4\le a\le3k$, $a\equiv 2$  mod 3 and $b= n - 2a + 4,$}  \\
    a+b-1  &   \text{if $4\le a\le3k$, $a\equiv 2$  mod 3 and $b=n - 2a +5,$ }  \\
    a+b & \text{otherwise.}
  \end{cases}
\end{equation}

\begin{lemma}\label{thirdlemma}
The latin square $\VV_n$, where $n=6k+4$, has a transversal, and all
suitable diagonals include the $k$ entries
\begin{align}\label{4modd6}
(1,0,3), (5,n-6,0), (8,n-12,n-3), (11,n-18,n-6),\ldots,( 3k-1,10,3k+10).
\end{align}
\end{lemma}

\begin{proof}
  Inspecting \cref{ThirdFamily}, we have $-3\le\Delta(r,c,s)\le3$ for every entry $(r,c,s)$ in $\VV_n$.
The sum of the minimum $\Delta$-value from each row is
\begin{align*}
0-1+0-3+\sum_{i=1}^{k-1} (0-1-2)= -3k-1 > -3k - 2 = - \dfrac{n}{2},
\end{align*}
whereas the sum of the maximum $\Delta$-value from each row is
$$3+2+0+0+\sum_{i=1}^{k-1} (2+1+0)=3k+2 = \dfrac{n}{2}.$$
Therefore, any suitable diagonal must include an entry with the maximum $\Delta$-value from each row, which implies that it includes the entries in \cref{4modd6}.

We next show that there exists a transversal $T$ in $\VV_n$, defined by
\begin{equation*}
\col(a)= \begin{cases}
  2  &   \text{if $a=0,$} \\
  0& \text{if $a = 1,$} \\
  6& \text{if $a = 2,$} \\
  1& \text{if $a = 3,$}\\
  4& \text{if $a = 3k+5,$} \\
  n-2a+5& \text{if $a \in \{ 3k+1,3k+2,3k+3\},$} \\
  n-2a+4   &    \text{if $4 \leq a \leq 3k$ and $a\not\equiv 0$  mod 3,}  \\  
  n-2a+7  &    \text{if $4 \leq a \leq 3k$ and $a\equiv 0$ mod 3,} \\  
  n-2a+6   &    \text{if $3k+6\leq a <n$ and  $a\equiv 0$ mod 3,}\\  
  n-2a+3 &    \text{if $3k+4\leq a <n$ and $a\equiv 1$  mod 3,} \\    
  n-2a+9  &    \text{if $3k+8\leq a <n$ and  $a\equiv 2$  mod 3.} 
\end{cases}
\end{equation*}
In order to prove that $T$ is a transversal, we show that every column
and symbol is represented exactly once in $T$.
Let $C_{a,b}$ and $S_{a,b}$ be, respectively, the set of columns and symbols used by the entries in $T$ in rows $\{a,a+1,\ldots,b\}$. We have
\begin{align*}
  C_{0,3}&=\{0,1,2,6\},\\
  C_{4,3k} &= \{6s+r: 1 \leq s \leq k - 1 , 4\le r\le6\},\\
  C_{3k+1,3k+3}&=\{3,5,7\},\\
  C_{3k+4,6k+3} &= \{6s+r: 2 \leq s \leq k  , 1\leq r\leq 3 \} \cup \{4,8,9\}.
\end{align*}
And 
\begin{align*}
 S_{0,3}&=\{3,4,5,8\},\\
S_{4,3k+3}&=\{s:3k+6\le s\le6k+6\} \setminus \{3k+9\},\\
S_{3k+4,n-1}&= \{s:9\le s\le3k+5\} \cup \{ 6,7,3k+9\}. 
\end{align*}
It is routine to check that $C_{0,n-1} = S_{0,n-1}=\{0,1,\ldots,n-1\}$,
completing the proof. 
\end{proof}

Combining \Cref{zeromodsix,twomodsix,thirdlemma} proves
\Cref{t:surrogate}.  We conclude this section by showing that our
latin squares constructed in the proof of \Cref{t:surrogate} contain
many transversal-free entries.

\begin{lemma}\label{initialcor}
For all $n\equiv 0 \bmod 6$, there exists a latin square $L$ of
order $n$ with $n^2>\tau(L) \geq (19n^2-51n+36)/36$.
\end{lemma}

\begin{proof}
If $n=6$ then we can take $L$ to be 
\begin{align*}
L=\begin{bmatrix}
\myredshadedcell{0} & 1 &2 &\myblueshadedcell{3} &4 &\myblueshadedcell{5} \\
1 &\myblueshadedcell{0}  &\myredshadedcell{3} &4 &\myblueshadedcell{5} &\myblueshadedcell{2} \\
2  &3 &\myblueshadedcell{1} &\myblueshadedcell{5} &\myblueshadedcell{0} &\myredshadedcell{4} \\
\myblueshadedcell{3} &\myredshadedcell{5} &\myblueshadedcell{4} &\myblueshadedcell{1} &2 &0 \\
\myblueshadedcell{4} &\myblueshadedcell{2} &5 &0 &\myredshadedcell{1} &\myblueshadedcell{3} \\
5 &\myblueshadedcell{4} &\myblueshadedcell{0} &\myredshadedcell{2} &3 &1  
\end{bmatrix}.
\end{align*}
Here the $16$ blue entries are transversal-free and the red entries
form a transversal.
Let $n=6k\geq 12$. Let $F$ be the set of $k$ pinned entries of $\TT_n$
given in \cref{0modd6}, and
\begin{align}\label{TSets}
C=\{6u:2\leq u \leq k\} \cup \{1\}, S=\{3u:k+3\leq u \leq 2k+1\} \cup \{4\},
\end{align}
to be, respectively, the set of columns and symbols of entries in $F$. Define 
\begin{align*}
N &=\{(r,c,s)\in \TT_n : c \in C \} \setminus F, \\
O &=\{(r,c,s)\in \TT_n : s \in S \} \setminus F.
\end{align*} 
Note that $s\equiv 0\, \text{or} \, 1\bmod 3$ for each entry $e=(r,c,s)\in O$. We also introduce the following sets:
\begin{align*}
P&=\{(r,c,s) \in \TT_n : r = 3u-2, 1\leq  u \leq k-1, \Delta(r,c,s)\neq2 \},\\
Q &=\{(r,c,s) \in \TT_n : r = 3u-1, 1\leq  u \leq k-1, \Delta(r,c,s)\neq1\}, \\
R &=\{(r,c,s) \in \TT_n : r = 3u, 2\leq  u \leq k-1, \Delta(r,c,s)\neq0\}.
\end{align*}
The proof of \Cref{zeromodsix} shows that the entries in $M=P\cup Q\cup R$ are transversal-free as they do not have the maximum $\Delta$-value in their row. 
Note that the entries in $N \cup O$ are also transversal-free.  We will calculate $|M\cup N \cup O|$ using inclusion-exclusion. To begin with, \Cref{zeromodsix} and \Cref{Structure}
imply that $|N| = |O| = k (n-1)$ and $|M|=|P|+|Q|+|R|=n-1+(k-2)(n/2+1)+(k-1)(n-1)+(k-2)n/2=2nk-2n-2$. Let us next examine the pairwise intersections, which fall naturally into three cases.

\begin{enumerate}[label = \textbf{Case \arabic*: }, ref =
    \textbf{Case \arabic*}, leftmargin=0pt, labelsep=0pt, itemindent=*]

\item \label{case1} 
  $M\cap N$.
  
  We evaluate how many entries of $N$ lie in $P$, $R$, and $Q$ respectively. Note that row $r=1$ shares $k-1$ entries with $N$. Suppose that $r>1$.  The structure of $\TT_n$ in \cref{Structure} illustrates that $c$ is odd for each entry $e=(r,c,s)\in P$, except the entries
\begin{equation}\label{e:excentT}
  \{(3u-2,6(k-u+1),n-3u+5):2\le u\le k-1\},
\end{equation}
for which $s\equiv r+c+1 \bmod n$. So the entries \eref{e:excentT}, as well as those of the form $(3u-2,1,3u-1)$ are in $P\cap N$. As a result, $|P\cap N|=3k-5$. In a similar manner we can see that $c$ is even for each entry $e=(r,c,s) \in R$, and that $|R\cap N|=(k-1)(k-2)$. Also, $|Q\cap N|=(k-1)^2$ follows from the definitions of $Q$ and $N$. We thus get $|M\cap N|=2(k^2-k-1)$.

\vspace*{0.1cm}
\item \label{case2}
  $M\cap O$.
  
  We evaluate how many entries of $O$ lie in $P$, $R$, and $Q$ respectively. It is straightforward to see that row $r=1$ shares $k-1$ entries with $O$. Assume $r>1$. The structure of $\TT_n$ in \Cref{Structure} reveals that $s\equiv r+c \bmod n$ and $c$ is odd for each entry $e\in P$ except for \eref{e:excentT} and
\begin{equation}\label{e:othexcentT}
  \{(3u-2,6(k-u+1)+1,n-3u+6):2\le u\le k-1\}.
\end{equation}
Note that $O$ contains \eref{e:othexcentT} but is disjoint from \eref{e:excentT}.
We also have $s\equiv r+c-2 \bmod n$ where $c$ is even for each entry
$e=(r,c,s) \in R$. Hence, with two exceptions arising from
\eref{e:excentT} and \eref{e:othexcentT}, every symbol occurs
in $P\cap O$ in row $r$ if and only if it does not appear in $R\cap O$
in row $r+2$. The two exceptions are symbol $n-3r\in O$, which occurs in
both places and $n-3r-1\notin O$ which occurs in neither.
Consequently, $|P\cap O| + |R\cap O|=k-1+(k-2)(k-1+2)$.  The definitions of $Q$ and $N$ indicate $|Q \cap O| = (k-1)^2$, leading to $|M\cap O|=2(k^2-k-1)$.

\vspace*{0.1cm}
\item \label{case3}
  $N\cap O$.
  
Given the latin property of square $\TT_n$ along with the definitions
of $N$ and $O$, one can see that $|N\cap O|=k(k-1)$.  We now look at
how many of these $k(k-1)$ entries lie in $P$, $Q$, and $R$ respectively.

 \begin{enumerate}[label=\textbf{Subcase \arabic{enumi}.\alph*: },ref=\textbf{Subcase \arabic{enumi}.\alph*},leftmargin=0pt,labelsep=0pt,itemindent=73pt,listparindent=\parindent]
\item\label{subcase1a} 
 $P\cap N \cap O$. 

  Let $e=(r,c,s)\in P\cap N \cap O$. First consider $r=1$. Note that 
  $\{(1,0,3),(1,2,2),(1,3,4)\}\cap N=\emptyset$ and $(1,1,1)\notin O$.
  All other entries with $r=1$ satisfy $s\equiv 1+c \bmod n$ and hence cannot
  have $s\=c\=0\bmod3$.
  Thus, row $r=1$ contains no entry in $N \cap O$. Next, we assume that $r>1$. By the definition of $P$, we have $s\equiv r+c\equiv1+c\bmod 3$ except for entries
  \eref{e:excentT} and \eref{e:othexcentT}, which can be ignored since they are disjoint from $O$ and $N$ respectively.
  Inspecting \eref{TSets}, the only possibility is that $s=4$ and $c\ne1$. 
  However, \cref{Structure} shows that $c\equiv 3 \,\text{or}\, 4 \bmod 6$ for $(3u-2,c,4)$ implying that $(3u-2,c,4)\notin N$.
  Therefore, $|P\cap N \cap O|=0$.

\vspace*{0.1cm}
\item \label{subcase1b} 
  $Q\cap N \cap O$.
   
  Consider $e = (r, c, s) \in Q \cap N \cap O$. Then $e\in Q$ implies that
  $s \equiv r + c\=2+c\bmod{3}$ holds, except for the entry $(3u-1,6(k-u+1)+1,r+c-1)$, which is not in $N$. 
  Inspecting \eref{TSets}, we conclude that $c=1$.
  So consider an entry in $Q \cap N$ of the form $e=(3u-1,1,3u)$, where $1\leq u \leq k-1$. If $u=1$, then $(2,1,3)\notin \TT_n$. If $u\geq 2$, then
$3u$ belongs to the set $ \{6,9,\ldots,3k-3\}$, which shares no symbols with $S$, thus $e\notin O$. Therefore, $|Q\cap N \cap O|=0$.

\vspace*{0.1cm}
\item \label{subcase1c} 
  $R\cap N \cap O$.
  
Assume $e=(r,c,s)\in R\cap N \cap O$. Thus, $s\equiv r+c-2\equiv c-2\bmod 3$ and $c\equiv0 \bmod 6$, unless $c=1$. Given \cref{TSets} and the definition of $R$, the only entries that $R$ shares with $N \cap O$ are of the form $(r,c,4)$ whenever $r$ is even. Thus, $|R\cap N \cap O|=\big\lceil (k-2)/2 \big\rceil$.
 \end{enumerate}
\end{enumerate}
Applying inclusion-exclusion reveals that 
$|M\cup N \cup O|\ge(19n^2-51n+36)/36$,  proving the lemma.
\end{proof}

We have the following lemma for orders $n \equiv 4 \bmod 6$. The
approach will be similar to \Cref{initialcor}, but we will need to be
more careful with the mod 3 argument given that $n\not\=0\bmod6$.  For
the purpose of all mod 3 arguments, we will treat the index set as
consisting of the \emph{integers} $0,1\dots,n-1$. For other arguments,
the index set will be taken modulo $n$.

\begin{lemma}\label{seccor}
Let  $n\equiv 4 \bmod 6$, where $n \geq 10$. There exists a latin square $L$ of order $n$ with $n^2>\tau(L) \geq  (19n^2-86n-68)/36$.
\end{lemma}

\begin{proof}
  Let $k= \big\lfloor {n}/{6} \big\rfloor$. 
Let $F$ be the set of $k$ pinned entries of $\VV_n$ given in \cref{4modd6}, and 
\begin{align}\label{USets}
C=\{6(k-u+1)+4:1\leq u \leq k\}, S=\{3u+1:k+3\leq u \leq 2k+1\} \cup \{3\},
\end{align}
be, respectively, the set of columns and symbols of entries in $F$. Define 
\begin{align*}
N &=\{(r,c,s)\in \VV_n : c \in C \} \setminus F, \\
O &=\{(r,c,s)\in \VV_n : s \in S \} \setminus F.
\end{align*}
Note that $s\equiv 1\bmod 3$ or $s=3$ for each entry $e=(r,c,s)\in O$.  We also introduce the following sets:
\begin{align*}
P&=\{(r,c,s) \in \VV_n : r = 3u-2, 1\leq  u \leq k, \Delta(r,c,s)\neq2 \},\\
Q &=\{(r,c,s) \in \VV_n : r = 3u-1, 2\leq  u \leq k, \Delta(r,c,s)\neq1\}, \\
R &=\{(r,c,s) \in \VV_n : r = 3u, 2\leq  u \leq k, \Delta(r,c,s)\neq0\}.
\end{align*}
The proof of \Cref{thirdlemma} shows that the entries in $M=P\cup Q\cup R$ are transversal-free as they do not contain the maximum $\Delta$-value in the row where they are located. 
Additionally, the entries in $N \cup O$ are transversal-free.  We now calculate $|M\cup N \cup O|$ using inclusion-exclusion. To start with, \Cref{thirdlemma} and \Cref{ThirdFamily}
imply that $|N| = |O| = k (n-1)$ and $|M|=|P|+|Q|+|R|=n-1+(k-1)(n/2+1)+(k-1)(n-1)+(k-1)n/2=2kn-n-1$. Next, we study the pairwise intersections in three cases.

\begin{enumerate}[label = \textbf{Case \arabic*: }, ref =
    \textbf{Case \arabic*}, leftmargin=0pt, labelsep=0pt, itemindent=*]
\item \label{C1} 
$M\cap N$.

  We evaluate how many entries of $N$ lie in $P$, $R$, and $Q$ respectively. Note that row $r=1$ shares $k-1$ entries with $N$. Suppose that $r>1$. The structure of $\VV_n$ in \cref{ThirdFamily} illustrates that $c$ is odd for each entry $e=(r,c,s)\in P$ except the entries
\begin{equation}\label{e:excentV}
  \{(3u-2,6(k-u+1)+4,n-3u+5):2\le u\le k\}
\end{equation}
for which $s\equiv r+c+1 \bmod n$. The entries in \eref{e:excentV} belong to $P\cap N$, and as a result, $|P\cap N|=2k-2$. In a similar manner we can see that $c$ is even  for each entry $e=(r,c,s) \in R$, and that $|R\cap N|=k(k-1)$. Also, $|Q\cap N|=(k-1)^2$ follows from the definitions of $Q$ and $N$. We thus get $|M\cap N|=2k^2-k-1$.

\vspace*{0.1cm}
\item \label{C2}
  $M\cap O$.
  
  We evaluate how many entries of $O$ lie in $P$, $R$, and $Q$ respectively. One can check that the row $r=1$ shares $k-1$ entries with $O$. Assume $r>1$. The structure of $\VV_n$ in \Cref{ThirdFamily} reveals that $s\equiv r+c \bmod n$ and $c$ is odd for each entry $e\in P$ except for those in
\eref{e:excentV} and
\begin{equation}\label{e:othexcentV}
  \{(3u-2,6(k-u+1)+5,n-3u+6):2\le u\le k\}.
\end{equation}
Note that $O$ contains \eref{e:othexcentV} but is disjoint from
\eref{e:excentV}.
We also have $s\equiv r+c-2 \bmod n$ where $c$ is even for each entry
$e=(r,c,s) \in R$. Again, each symbol in $S\setminus\{n-3r\}$ occurs in
$P\cap O$ in
row $r$ if and only if it does not appear in $R\cap O$ in row $r+2$.
The exceptional symbol $n-3r$ gets counted in both rows.
As a result, $|P\cap O| + |R\cap O|=(k-1)(1+k-1+2)$.
Also, the definitions of $Q$ and $N$ indicate $|Q\cap O|=(k-1)^2$.
Consequently, $|M\cap O|=2k^2-k-1$. 

\vspace*{0.1cm}
\item \label{C3}
  $N\cap O$.
  
Given the latin property of square $\VV_n$ along with the definitions
of $N$ and $O$, one can see that $|N\cap O|=k(k-1)$. We now look at how
many of these $k(k-1)$ entries lie in $P$, $Q$, and $R$ respectively.

\begin{enumerate}[label=\textbf{Subcase \arabic{enumi}.\alph*: }, ref = \textbf{Subcase \arabic{enumi}.\alph*}, leftmargin=0pt, labelsep=0pt, itemindent=73pt, listparindent=\parindent]
\item \label{sbcse1a}
  $P\cap N \cap O$. 
  
Let $e=(r,c,s)\in P\cap N \cap O$. If $r=1$, then \cref{ThirdFamily}
shows that $s\equiv r+c \bmod n$ except where $c\in\{0,1,2\}$, but
$(1,0,3)\in F$ and $C\cap \{1,2\} =\varnothing$.
Thus we may assume that $s\equiv 1+c \equiv2 \bmod 3$,
which means $e\notin O$.
Thus, row $r=1$ shares no entries with
$P\cap N \cap O$. Assume $r>1$.
The entries in \eref{e:excentV} are not in $O$ because
$n-3u+5 \notin S$. All other entries in $P$ have odd $c\notin C$.
Therefore, $|P\cap N \cap O|=0$.

\vspace*{0.1cm}
 \item \label{sbcse1b} 
   $Q\cap N \cap O$.
   
Let $e=(r,c,s)\in Q\cap N \cap O$.  Then $e\in Q\cap N$ yields that $s\equiv r+c \bmod n$. First we show the column $c=0$ shares no entries with $Q\cap N \cap O$. If $c=0$, then $e$ must belong to the set $\{(3u-1,0,3u-1):2\leq u\leq k\}$, which is disjoint from $O$ because $5\leq s \leq 3k-1$.
  Next assume $c\neq 0$. Hence $s\equiv r+c \bmod n$ implies that in $\Z$, either $s=r+c$ or $s+n=r+c$. For the former case we have $s\equiv 0\bmod 3$ and $s\geq 5$, so $e\notin O$. In the latter case $s\equiv 2\bmod3$ which implies  $e\notin O$. Therefore, $|Q\cap N \cap O|=0$.

\vspace*{0.1cm}
\item \label{sbcse1c} 
  $R\cap N \cap O$.
  
Assume $e=(r,c,s)\in R\cap N \cap O$. Thus $e\in R$ yields that $s\equiv r+c-2 \bmod n$. First we examine the case where $c=0$. In this case, $e$ must belong to the set 
$\{(3u,0,3u-2): 2\leq u \leq k\}$
 which is disjoint from $O$ because $4\le s\le3k-2$. Thus $e\notin O$.  
Now assume $c\neq 0$. 
 Note that $s\equiv r+c-2 \bmod n$ implies that in $\Z$ either $s= r+c-2$ or $s+n= r+c-2$. We will show that $|R\cap N \cap O|=0$ by analysing both options. 
If $s=r+c-2\equiv c+1\bmod3$, then $s \equiv 2 \bmod3$, so $e\notin O$. 
Next assume $s=r+c-2-n$. Hence $s\equiv c\equiv 1\bmod 3$
and $s\le3k+n-6-2-n=3k-8$, 
so $e\notin O$. Therefore, $|R\cap N \cap O|=0$.
\end{enumerate}

\end{enumerate}
Applying inclusion-exclusion reveals that
$|M\cup N \cup O|=(19n^2-86n-68)/36$, and the proof of the lemma is complete.
\end{proof}

Lastly, the following lemma covers all remaining orders that are
$n\equiv 2 \bmod 6$.

\begin{lemma}\label{thirdcor}
Let $n\equiv 2 \bmod 6$, where $n\geq 8$. There exists a latin square $L$ of order $n$ with $n^2>\tau(L) \geq  (19n^2-73n-182)/36$.
\end{lemma}

\begin{proof}
If $n=8$ then we can take 
\begin{align*}
L=\begin{bmatrix}
0&1&2&3&\myredshadedcell{4}&5&6&7\\
1&0&3&2&5&4&7&\myredshadedcell{6}\\
2&\myredshadedcell{3}&\myblueshadedcell{1}&\myblueshadedcell{0}&\myblueshadedcell{7}&\myblueshadedcell{6}&4&5\\
3&2&\myblueshadedcell{6}&\myblueshadedcell{7}&\myblueshadedcell{0}&\myblueshadedcell{1}&\myredshadedcell{5}&4\\
\myblueshadedcell{4}&\myblueshadedcell{7}&\myblueshadedcell{5}&\myredshadedcell{1}&6&\myblueshadedcell{3}&2&\myblueshadedcell{0}\\
\myblueshadedcell{5}&\myblueshadedcell{4}&\myblueshadedcell{7}&\myblueshadedcell{6}&3&\myredshadedcell{2}&0&1\\
\myblueshadedcell{6}&\myblueshadedcell{5}&\myredshadedcell{0}&\myblueshadedcell{4}&\myblueshadedcell{2}&7&\myblueshadedcell{1}&\myblueshadedcell{3}\\
\myredshadedcell{7}&\myblueshadedcell{6}&\myblueshadedcell{4}&5&1&0&3&2
\end{bmatrix},
\end{align*}
where the $25$ blue entries are transversal-free and the red entries
form a transversal.

Let $k=\big\lfloor n/6 \big\rfloor$. We again utilise a similar
argument to \Cref{initialcor} and \Cref{seccor}.  Let $F$ be the set
of $k$ pinned entries of $\UU_n$ given in \Cref{twomodsix}, and
\begin{align*}
C=\{6u+2:2\leq u \leq k-1\} \cup \{3,4\}, S=\{3u:k+4 \leq u \leq 2k+1\} \cup \{5,8\},
\end{align*}
be, respectively, the set of columns and symbols of entries in $F$. Define 
\begin{align*}
N &=\{(r,c,s)\in \UU_n : c \in C \} \setminus F, \\
O &=\{(r,c,s)\in \UU_n : s \in S \} \setminus F.
\end{align*}
Note that $s\equiv 0\bmod 3$ or $s\in \{1,5,8\}$ for each entry $e=(r,c,s)\in O$. We also introduce the following sets:
\begin{align*}
P&=\{(r,c,s) \in \UU_n : r = 3u-1, 2\leq  u \leq k-1, \Delta(r,c,s)\neq2 \},\\ \nonumber
Q &=\{(r,c,s) \in \UU_n : r=1\text{ or }r = 3u\text{ for }1\leq u \leq k-1; \Delta(r,c,s)\neq1\}, \\ \nonumber
R &=\{(r,c,s) \in \UU_n : r = 3u+1, 2\leq  u \leq k-1, \Delta(r,c,s)\neq0\}. 
\end{align*}
The proof of \Cref{twomodsix} shows that entries in $M=P\cup Q\cup R$ are transversal-free as they do not contain the maximum $\Delta$-value in the row where they are located. 
Additionally, the entries in $N \cup O$ are also transversal-free.  We now calculate $|M\cup N \cup O|$ using inclusion-exclusion. To start with, \Cref{twomodsix} and \Cref{SecFamily}
imply that $|N| = |O| = k (n-1)$ and $|M|=|P|+|Q|+|R|=(k-2)(n/2+1)+k(n-1)+(k-2)n/2=2kn-2n-2$. Next, we study the pairwise intersections in three cases.
\begin{enumerate}[label = \textbf{Case \arabic*: }, ref =
    \textbf{Case \arabic*}, leftmargin=0pt, labelsep=0pt, itemindent=*]
\item \label{step1} 
  $M\cap N$.
  
We evaluate how many entries of $N$ lie in $P$, $R$, and $Q$
respectively. The structure of $\UU_n$ in \cref{SecFamily} illustrates
that $c$ is odd for each entry $e=(r,c,s)\in P$, except the entries
\begin{equation}\label{e:excentU}
\{(3u-1,6(k-u+1)+2,6k-3u+8):2\le u\le k-1\},
\end{equation}
for which $s=r+c+1$. The entries in \eref{e:excentU} belong to $P\cap N$, as do the entries of the form $(3u-1,3,3u+2)$. As a result, $|P\cap N|=2(k-2)$. In a similar manner, we can see that $c$ is even for each entry $e=(r,c,s) \in R$, and that $|R\cap N|=(k-1)(k-2)$. Also, $|Q\cap N|=k(k-1)$ follows from the definitions of $Q$ and $N$. We thus get $|M\cap N|=2k^2-2k-2$.

\vspace*{0.1cm}
\item \label{step2}
  $M\cap O$.
  
  We evaluate how many entries of $O$ lie in $P$, $R$, and $Q$ respectively. The structure of $\UU_n$ in \cref{SecFamily} reveals that $s\equiv r+c \bmod n$ and $c$ is odd for each entry $e\in P$ except for those in \eref{e:excentU} and
  \begin{equation}\label{e:othexcentU}
    \{(3u-1,6(k-u+1)+3,6k-3u+9):2\le u\le k-1\}.
  \end{equation}
Note that $O$ contains \eref{e:othexcentU} but is disjoint from
\eref{e:excentU}.
We also have $s\equiv r+c-2 \bmod n$ where $c$ is even for each entry
$e=(r,c,s) \in R$. Again, any symbol in $S\setminus\{n-3r+6\}$
occurs in $P\cap O$ in
row $r$ if and only if it does not appear in $R\cap O$ in row $r+2$.
The exceptional symbol $n-3r+6$ counts in both rows.
As a result, $|P\cap O| + |R\cap O|=(k-2)(k+1)$. Finally, the definitions of $Q$ and $N$ indicate $|Q\cap O|=k(k-1)$. In conclusion, $|M\cap O|=2k^2-2k-2$. 

\vspace*{0.1cm}
\item \label{step3}
  $N\cap O$.
  
Given the latin property of square $\UU_n$ along with the definitions
of $N$ and $O$, one can see that $|N\cap O|=k(k-1)$.  We now look at
how many of these $k(k-1)$ entries lie in $P$, $Q$, and $R$
respectively.

\begin{enumerate}[label = \textbf{Subcase \arabic{enumi}.\alph*: }, ref = \textbf{Subcase \arabic{enumi}.\alph*}, leftmargin=0pt, labelsep=0pt, itemindent=73pt, listparindent=\parindent]
\item \label{subc1a} 
  $P\cap N \cap O$.
  
  Let $e=(r,c,s)\in P\cap N \cap O$. Thus $e\in P$ implies $s\equiv r+c\bmod n$ and $c$ is odd except for the entries in
\eref{e:excentU} where $s\equiv r+c+1 \bmod n$. First we analyse the case $c=3$. In this case, $e$ belongs to the set 
$\{(3u-1,3,3u+2): 2\leq u \leq k - 1\}$. Inspection shows that $(5,3,8)$ is the only entry of this set that belongs to $O$. Now we assume that $c\neq 3$. For each entry in \eref{e:excentU}, $s=6k-3u+8 \equiv 2\bmod 3$ and $3k+11\leq s\leq n$, so $e\notin O$. Therefore $|P\cap N \cap O|=1$.

  \vspace*{0.1cm}
 \item \label{subc1b} 
   $Q\cap N \cap O$.
   
Consider $e = (r, c, s) \in Q \cap N \cap O$. Thus $e\in Q\cap N$ reveals that $s\equiv r+c \bmod n$. We analyse cases $r=1$, $r=3$, and $r\notin\{1,3\}$ separately. First assume $r=1$. If $c\in \{3,4\}$, then $e\in \{(1,3,5),(1,4,4)\}$ whereas $(1,3,5) \in F$ and $(1,4,4)\notin O$. Therefore we assume $c\notin  \{3,4\}$. In this case, the $\big\lfloor (k-3)/2 \big\rfloor$ entries in
$\{(1,6u+2,6u+3): 2\leq u \leq k-1\} \cap S$ 
belong to $N\cap O$.

Next assume $r=3$. If $c\in \{3,4\}$, then $e\in \{(3,3,4),(3,4,8)\}$ whilst $(3,3,4)\notin O$ and $(3,4,8)\in F$. Therefore we can assume $c\notin \{3,4\}$. In this case $s \equiv 3+c \equiv 2\bmod 3$ and $17 \leq s\leq n - 3$, so $e\notin O$.

Now we study the case where $r\notin\{1,3\}$. If  $c\in \{3,4\}$, then $9 \leq s \leq 3k+1$ which is outside the range of symbols in $S$.
Hence we assume $c\notin \{3,4\}$. 
Note that $s\equiv r+c \bmod n$ implies in $\Z$ either $s = r + c $ or
$s+n = r + c$. In the former case, $s\equiv 2\bmod 3$ and $s\geq20$,
so $e\notin O$. In the latter case we notice that $s\equiv 0 \bmod 3$
and $s\leq3k-3+6k-4-n=3k-9$, so $e\notin O$. As a result
$|Q\cap N \cap O|=\big\lfloor (k-3)/2 \big\rfloor\ge(k-4)/2$.

\vspace*{0.1cm}
\item \label{subc1c} 
  $R\cap N \cap O$.
  
Assume $e=(r,c,s)\in R\cap N \cap O$. Hence, $e\in R$ implies that $s\equiv r+c-2 \bmod n$ and $c$ is even. If $c=4$ then $e$ must belong to the set
$\{(3u+1,4,3u+3): 2\leq u \leq k-1\},$
which is disjoint from $O$. So we can assume that $c\neq4$. The
condition $s\equiv r+c-2 \bmod n$ implies that in $\Z$ either
$s=r+c-2$ or $s+n=r+c-2$. For the first case we have $8\le s\equiv1\bmod3$,
so $e\notin O$. In the latter case $s\equiv 2\bmod3$, so
the only possible values for $s$ are $5$ and $8$.
For odd $r>12$ there is a $c$ such that $(r,c,5)\in R\cap N$  and
for even $r>12$ there is a $c$ such that $(r,c,8)\in R\cap N$.
 Thus $|R\cap N \cap O|=k-4$.
\end{enumerate}
\end{enumerate} 
Applying inclusion-exclusion reveals that
$|M\cup N \cup O|\geq (19n^2-73n-182)/36$, and the proof of the lemma is complete.
\end{proof}

Combining \Cref{initialcor,seccor,thirdcor}, we have covered all cases
needed to prove \Cref{resultfrommain}.  The families $\TT_n$, $\UU_n$,
and $\VV_n$ each achieve at least $19n^2/36+O(n)$ transversal-free
entries.  For $n\geq 88$, more than half of the entries are
transversal-free.  We have not worked hard to count the number of
transversal-free entries precisely, in part because it would be
tedious to show that all other entries are in transversals.
Nevertheless, we believe that the bounds we have given are within
$O(n)$ of the true number of transversal-free entries. Some evidence
for this belief is provided by the computations reported in
\Tref{T:howclose}. The table shows our lower bound on the number, and
the true number of transversal-free entries for each order in the
range $10\le n\le24$.

\begin{table}
\begin{tabular}{|c|c|c|}
  \hline
  $L$& Lower bound on $\tau(L)$& Actual $\tau(L)$\\
  \hline
$V_{10}$& 27& 34\\
$T_{12}$& 60& 67\\
$U_{14}$& 70& 88\\
$V_{16}$& 95& 107\\
$T_{18}$& 147& 159\\
$U_{20}$& 166& 190\\
$V_{22}$& 201& 217\\
$T_{24}$& 271& 287\\
\hline
\end{tabular}
\medskip
\caption{\label{T:howclose}Number of transversal-free entries in our latin squares.}
\end{table}

\section{Latin squares with special subsquares}\label{secsec}

Available evidence suggests that it is much harder to restrict
transversals for latin squares of odd orders.  For orders $n\in
\{5,7\}$, there is a latin square of order $n$ with a pinned entry
(see, e.g.~\citep{egan2012latin}). There is no such latin square known
for any larger odd order, and for order 9 they definitely do not
exist.  In this section we construct latin squares of odd order where
every transversal has to hit specific subsquares rather than specific
entries.

Let $L_1$ and $L_2$ be two latin squares of order $n$. An ordered
triple $\phi=(\alpha,\beta,\gamma)$ of permutations $\alpha$, $\beta$,
and $\gamma$ is called an \emph{isotopism} from $L_1$ to $L_2$ if $L_2$
is obtained by permuting rows, columns and symbols of $L_1$ by
$\alpha$, $\beta$, and $\gamma$, respectively. If $L_1=L_2$ then
$\phi$ is an \emph{autotopism}. An autotopism with the property
that $\alpha=\beta=\gamma$ is called an \emph{automorphism}.

For odd $m \geq 3$, we define the latin square $\LL_n$ of order $n=3m$
and indexed by $\Z_n$ by 
\begin{align*}
\LL_n=\begin{bmatrix}
	\AA_{11}&\AA_{12}&\AA_{13}\\
	\AA_{21}&\AA_{22}&\AA_{23}\\
	\AA_{31}&\AA_{32}&\AA_{33}\\
\end{bmatrix},
\end{align*}
where each block $\AA_{ij}$ is an $m\times m$ latin subsquare
defined by 
\begin{align*}
  \AA_{ij}[a,b] = \begin{cases}
    a+b    &    \text{if $ a+b \not\equiv 0$  mod $m$  and  $i+j\equiv 2$  mod 3,}\\
    a+b+m    &    \text{if $a+b \not\equiv m -1$  mod $m$  and  $i+j\equiv0$  mod 3,} \\
    a+b+2m    &    \text{if $a+b \not\equiv m-1$  mod $m$  and $i+j\equiv1$  mod 3,} \\
    0 &   \text{if $a+b\equiv0$  mod $m$  and  $(i,j)=(1,1),$} \\
    2m-1 &  \text{if $a+b\equiv0$  mod $m$  and  $(i,j)=(2,3),$}\\
    3m-1 & \text{if $a+b\equiv0 $  mod $m$  and $(i,j)=(3,2),$}  \\ 
    0 &   \text{if $a+b\equiv m-1$  mod $m$  and  $(i,j)\in \{(2,2),(3,3)\},$}\\
    2m-1 &  \text{if $a+b\equiv m-1$  mod $m$  and  $(i,j)\in\{(1,2),(3,1)\},$}\\
    3m-1 & \text{if $a+b\equiv m-1$ mod $m$  and  $(i,j)\in \{(1,3),(2,1)\}.$} 
  \end{cases}
\end{align*}

First, we introduce a function $\Delta_m : \LL_n \rightarrow \Z_m$
on the entries $e = (r,c,s)$ of $\LL_n$ by
$\Delta_m(r,c,s) \equiv r+c \bmod m$.
Assume that $T$ is a transversal in $\LL_n$. Then, by a similar
argument to \Cref{l:Delta}, we have
 \begin{align}\label{eq6}
   \sum_{(r,c,s)\in T} \Delta_m(r,c,s)=\sum_{r,c}(r+c)=2\sum_{i=0}^{n-1}i=n(n-1)
   \equiv 0 \bmod m.
 \end{align}
We define $x_{ij}$ to be the number of entries that $T$ includes from the subsquare $\AA_{ij}$, where $i,j\in\{1,2,3\}$. 
The need for $T$ to include one symbol from each row and column of $\LL_n$ implies 
\begin{align}
x_{11}+x_{12}+x_{13}  &=m, \label{e:row1}\\ 
x_{21}+x_{22}+x_{23}  &=m, \label{e:row2}\\
x_{31}+x_{32}+x_{33}  &=m, \label{e:row3}\\
x_{11}+x_{21}+x_{31}  &=m, \label{e:col1}\\ 
x_{12}+x_{22}+x_{32}  &=m, \label{e:col2}\\ 
x_{13}+x_{23}+x_{33}  &=m. \label{e:col3}
\end{align}
The value of $\Delta_m(e)$ for each $e=(r,c,s)$ of $\AA_{ij}$, is as follows:
\begin{align}\label{deltavalues}
 C=\begin{bmatrix}
 0&1&\ldots&m-2&m-1\\
 1&2&\ldots&m-1&0\\
 \vdots&\vdots&\ddots&\vdots&\vdots\\
 m-1&0&\ldots&m-3&m-2\\
 \end{bmatrix}.
\end{align}
Additionally, no matter from which block the entries of $T$ are
chosen, the sum of the $\Delta_m$-values of the entries containing
symbols $\{1,2,\ldots,m-1\}, \{m,m+1,\ldots,2m-2\}$ and
$\{2m,2m+2,\ldots,3m-2\}$ are, respectively, $1+2+\cdots+m-1$,
$0+1+\cdots+m-2$, and $0+1+\cdots+m-2$ using
(\ref{deltavalues}). Therefore, the sum of $\Delta_m$-values of all
entries in $T$ except the ones including the symbols in
$\{0,2m-1,3m-1\}$ is
\begin{align*}
(m-1)(m-2)+\dfrac{m(m-1)}{2} \equiv 2 \bmod  m.
\end{align*}
On the other hand, the $\Delta_m$-value of any entry including a
symbol in $\{0,2m-1,3m-1\}$ is either $0$ or $m-1$.  To satisfy
(\ref{eq6}), the entries of $T$ that include symbols in
$\{0,2m-1,3m-1\}$ must include one entry with $\Delta_m$-value of $0$
and two entries with $\Delta_m$-value of $m-1$. We will show that this
property forces all transversals of $\LL_n$ to contain at least one
entry from each subsquares $\AA_{ij}$ where $1\leq i,j\leq 3$.

\begin{lemma}\label{FirstSubsquare}
 Every transversal of $\LL_n$ must hit subsquare $\AA_{22}$ at least once.
\end{lemma}

\begin{proof}
  For the sake of contradiction, suppose that there exists a transversal $T$
  with $x_{22}=0$.
 From \eref{e:col1},
 \eref{e:col3} and \eref{e:row2} we then deduce that
 $x_{11}+x_{13}+x_{31}+x_{33}=m$. However, $T$ needs to include the
 $m$ symbols $\{0,2m,2m+1,\dots,3m-2\}$, which are now only available for
 selection in $\AA_{11}$, $\AA_{13}$, $\AA_{31}$ and $\AA_{33}$.  It
 follows that no other symbols can be chosen from these four blocks,
 and that $x_{11}+x_{33}=1$ and $x_{13}+x_{31}=m-1$.
 Now \eref{e:row3} and \eref{e:col3} imply that
 \begin{equation}\label{e:rc33}
   x_{23}+x_{32}=2m-x_{13}-x_{31}-2x_{33}=m+1-2x_{33}.
 \end{equation}
 Note that the $m-1$ symbols in $\{1,\dots,m-1\}$ must be chosen from
 $\AA_{23}\cup\AA_{32}$ and the only other symbols that are available
 in these blocks are $2m-1$ and $3m-1$. It follows that there are only
 two ways to satisfy \eref{e:rc33}. Either $x_{33}=0$ and we choose
 both of $2m-1$ and $3m-1$ from $\AA_{23}\cup\AA_{32}$, or $x_{33}=1$
 and we choose neither of $2m-1$ and $3m-1$ from
 $\AA_{23}\cup\AA_{32}$. Both of these options are incompatible with the
 requirements on the $\Delta_m$ values of the entries containing the
 symbols in $\{0,2m-1,3m-1\}$.
\end{proof}

Next, we show that every transversal has to contain at least one entry
from subsquare $\AA_{11}$.

\begin{lemma}\label{SecSubsquare}
Every transversal in $\LL_n$ contains at least one entry from the subsquare $\AA_{11}$.
\end{lemma}

\begin{proof}
  For the sake of contradiction, suppose that there exists a
  transversal $T$ for which $x_{11}=0$. From \eref{e:col2},
  \eref{e:col3} and \eref{e:row1} we then deduce that
  $x_{22}+x_{23}+x_{32}+x_{33}=m$. However, $T$ needs to include the
  $m$ symbols $\{0,1,\dots,m-1\}$, which are now only available for
  selection in $\AA_{22}$, $\AA_{23}$, $\AA_{32}$ and $\AA_{33}$.  It
  follows that no other symbols can be chosen from these four blocks.
  This forces all three symbols in $\{0,2m-1,3m-1\}$ to be chosen
  in entries where the $\Delta_m$ value is $m-1$, giving a
  contradiction.
\end{proof}

Now we prove \Cref{mySecTheorem}.

\begin{theorem}
Every transversal in $\LL_n$ contains at least one entry from each of
the subsquares $\AA_{11}, \AA_{12}, \AA_{13}, \AA_{21}, \AA_{22},
\AA_{23}, \AA_{31}, \AA_{32}$ and $\AA_{33}$.
\end{theorem}

\begin{proof}
It can be easily checked that $\tau=(\alpha,\alpha,\alpha)$ where
$\alpha$ is defined by
$$\alpha=(m, 2m)(m+1 ,2m+1)\cdots(2m-1,3m-1),$$ 
is an automorphism of $\LL_n$ with the following properties
\begin{align}\label{Aut1}
 \tau(\AA_{22})=\AA_{33}, \,\, \tau(\AA_{12})=\AA_{13},  \,\, \tau(\AA_{31})=\AA_{21},  \,\,  \tau(\AA_{23})=\AA_{32}.
 \end{align}
Thus, \Cref{FirstSubsquare} implies every transversal of $\LL_n$ hits
the subsquare $\AA_{33}$ at least once.

Also, $\phi=(\alpha,\beta,\gamma)$ where $\alpha$, $\beta$ and
$\gamma$ are defined as
\begin{align*}
\alpha &=(0,m)(1,m+1)(2,m+2)\cdots(m-1,2m-1), \\ 
\beta  &=(0,2m)(1,2m+1)(2,2m+2)\cdots(m-1,3m-1),  \\ 
\gamma  &=(0,2m-1)(m,2m)(m+1,2m+1)\cdots(2m-2,3m-2), 
\end{align*}
is an autotopism with the following properties 
\begin{align}\label{Auto2}
\phi(A_{11})=A_{23}, \,\,  \phi(A_{22})=A_{12}, \,\, \phi(A_{33})=A_{31}.
 \end{align}
Considering (\ref{Aut1}), (\ref{Auto2}), \Cref{FirstSubsquare}, and  \Cref{SecSubsquare}, it is straightforward to show that every transversal has to hit each of the subsquares $\AA_{11}, \AA_{12}, \AA_{13}, \AA_{21}, \AA_{22}, \AA_{23}, \AA_{31}, \AA_{32}$ and 
$\AA_{33}$ at least once.
\end{proof}

We finish the section by showing that $\LL_n$ has a transversal $T$.
If $n = 9$, then we may take $T$ to be
$$\{(0,1,1)(1,4,5)(2,7,6)(3,8,2)(4,0,4)(5,3,0)(6,6,3)(7,5,8)(8,2,7)\}.$$
For $n > 9$, we specify $T$ in two cases.
If $n \equiv 0 \bmod 9$ then we define $T$ by
\begin{equation*}
\col(a)= \begin{cases}
  2m/3 -2(a+1)  &   \text{if $0 \leq a < m/3,$}  \\
  4m/3 +a  &   \text{if $m/3 \leq a < 2m/3,$}  \\
  4m - 2a - 1  &   \text{if $2m/3 \leq  a < m,$}  \\
   a -  m/3   &   \text{if $m \leq a < 4m/3,$} \\
 13m/3 - 2a - 1 & \text{if $ 4m/3 \leq a < 5m/3,$} \\
 6m-2(a+1)& \text{if $ 5m/3 \leq a < 2m,$} \\
 14m/3 - 2a - 1  &   \text{if $2m \leq a < 7m/3,$} \\
  19m/3 - 2(a+1)& \text{if $ 7m/3 \leq a < 8m/3,$} \\
   a& \text{if $ 8m/3 \leq a < n.$} 
\end{cases}
\end{equation*}
Meanwhile, if $n \not\equiv 0 \bmod 9$ then we define $T$ by
\begin{equation*}
\col(a)= \begin{cases}
  -a(m+1)/2 &   \text{if $a\equiv 0 \bmod 3,$} \\
  -2a& \text{if $a\equiv 1 \bmod 3,$} \\
  m-a(m+1)/2& \text{if $a \equiv 2 \bmod 3.$}  
\end{cases}
\end{equation*}

\section{Concluding remarks}

In this paper we have constructed latin squares of even order that
have $\big\lfloor n/6\big\rfloor$ pinned entries.  It remains open
whether there are squares with more pinned entries. It is also open
for odd orders $>9$ whether there are latin squares that contain any
pinned entry.

It is known from \cite{Bal90} that latin squares of even order have an
even number of transversals. It follows that such squares cannot have
a unique transversal, and hence by \cite{CW10}, they cannot have more
than $n-3$ pinned entries. It remains open to prove a less trivial
upper bound on the number of pinned entries. Also, it seems highly
unlikely, but can a latin square of odd order have only one
transversal?

In terms of number of entries that are not in any transversals, we
have shown for even orders that asymptotically more than half of the
entries may be in no transversal, even when there are transversals
through some entries.  \Cref{t:quadodd} is the only result to exhibit
latin squares of odd order in which more than a constant fraction of
their entries are transversal-free.  A corresponding result is not
known for orders that are 1 or 5 mod~$6$.

Finally, for which orders are there latin squares that have
transversals, do not have any pair of disjoint transversals and yet do
not have any pinned entry? From \cite{egan2012latin}, we see
that no such square exists for orders less than 10, but we conjecture
that they exist for all large even orders.

  \let\oldthebibliography=\thebibliography
  \let\endoldthebibliography=\endthebibliography
  \renewenvironment{thebibliography}[1]{%
    \begin{oldthebibliography}{#1}%
      \setlength{\parskip}{0.2ex}%
      \setlength{\itemsep}{0.2ex}%
  }%
  {%
    \end{oldthebibliography}%
  }

\end{document}